\documentclass[12pt,reqno]{amsart} 

\usepackage[all]{xy}
\usepackage{supertabular}

\def\theoremnumberstyle{plain} 
\def\proofname{Proof}
\usepackage{tom}

\newcommand{\tom}[1]{}

\allowdisplaybreaks[3]

\SelectTips{cm}{12}  
\CompileMatrices

\begin{document}

   \parindent0cm

   \title[Reducible Families]{Reducible Families of Curves with Ordinary
     Multiple Points on Surfaces in $\PC^3$}
   \author{Thomas Keilen}
   \address{Universit\"at Kaiserslautern\\
     Fachbereich Mathematik\\
     Erwin-Schr\"odinger-Stra\ss e\\
     D -- 67663 Kaiserslautern
     }
   \email{keilen@mathematik.uni-kl.de}
   \urladdr{http://www.mathematik.uni-kl.de/\textasciitilde keilen}
   \thanks{The results were obtained during a research visit at the
     universities of Siena and Torino supported by the EAGER-node of North Italy.}


   \date{September, 2003.}

     
   \begin{abstract}
     In \cite{Kei03}, \cite{Kei05} and \cite{Kei04} we gave
     numerical conditions which ensure that an equisingular family is
     irreducible respectively T-smooth. Combining results from
     \cite{GLS01} and an idea from \cite{CC99} we give in the
     present paper series of examples of 
     families of irreducible curves on surfaces in $\PC^3$ with only ordinary
     multiple points which are reducible and where at least one
     component does not have the expected dimension. The examples show
     that for families of curves  with ordinary
     multiple points the conditions for T-smoothness 
     in \cite{Kei04} have the right asymptotics.
   \end{abstract}

   \maketitle

     Throughout this article $\Sigma$ will denote a smooth projective surface
     in $\PC^3$ of degree $n\geq 2$, and $H$ will be a hyperplane
     section of $\Sigma$. For a positive integer $m$ we denote by
     $M_m$ the topological singularity type of an ordinary $m$-fold
     point, i.\ e.\ the singularity has $m$ smooth branches with
     pairwise different tangents. And for positive integers $d$ and
     $r$ we denote by  $V_{|dH|}^{irr}(r M_m)$ 
     the family of irreducible curves in the linear system $|dH|$
     with precisely $r$ singular points all of which are ordinary
     $m$-fold points. $V_{|dH|}^{irr}(r M_m)$ is called T-smooth if it
     is smooth of the expected dimension
     \begin{displaymath}
       \expdim\left(V_{|dH|}^{irr}(r M_m)\right)=\dim|dH|-r\cdot \frac{m^2+m-4}{2}.
     \end{displaymath}

     \begin{theorem}\label{thm:thm}
       For $m\geq 18n$  there is an integer $l_0=l_0(m,\Sigma)$ such
       that for all $l\geq l_0$ the family
       $V_{|dH|}^{irr}(rM_m)$ with $d=2lm+l$ and $r=4l^2n$
       has at least one T-smooth component and one component of
       higher dimension. 

       Moreover, the T-smooth component dominates $\Sym^r(\Sigma)$
       under the map 
       \begin{displaymath}
         V_{|dH|}^{irr}(rM_m)\longrightarrow \Sym^r(\Sigma):
         C\mapsto \Sing(C)
       \end{displaymath}
       sending a curve $C$ to its singular locus, and the fundamental
       group $\pi_1(\Sigma\setminus C)$ of the complement of any curve $C\in
       V_{|dH|}^{irr}(rM_m)$ is abelian. 
     \end{theorem}
     
     Before we prove the theorem let us compare the result with the
     conditions for T-smoothness in \cite{Kei04} and for
     irreducibility in \cite{Kei05}. 

     Here we have given examples of non-T-smooth families $V_{|dH|}(rM_m)$ where 
     \begin{displaymath}
       r\cdot m^2 \equiv 
       n\cdot d^2,
     \end{displaymath}
     if we neglect the terms of lower order in $ml$. If $n\geq 4$ and the Picard
     number of $\Sigma$ is one, then according to
     \cite{Kei04} Corollary 2.3 respectively Corollary 2.4 -- neglecting terms of lower order
     in $m$ and $d$ --
     \begin{displaymath}
       r\cdot m^2
       < 
       \frac{1}{2n-6}\cdot n\cdot d^2
     \end{displaymath}
     would be a sufficient condition for T-smoothness. Similarly, if $n=2$, then $\Sigma$ is
     isomorphic to $\PC^1\times \PC^1$ and we may apply \cite{Kei04}
     Theorem 2.5 to find that 
     \begin{displaymath}
       r\cdot m^2<\frac{1}{8}\cdot n\cdot d^2
     \end{displaymath}
     implies T-smoothness. Since the families fail to satisfy the
     conditions only by a constant factor we see that asymptotically
     in $d$, $m$ and $r$ the conditions for T-smoothness are proper. 

     For irreducibility the situation is not quite as good. The
     conditions in \cite{Kei05} Corollary 2.4 for irreducibility if
     $n\geq 4$ and the Picard number of $\Sigma$ is one is roughly 
     \begin{displaymath}
       r\cdot m^2<\frac{24}{n^2 m^2}\cdot n\cdot d^2,
     \end{displaymath}
     and similarly for $n=2$ \cite{Kei05} Theorem 2.6 it is
     \begin{displaymath}
       r\cdot m^2<\frac{1}{6 m^2}\cdot n\cdot d^2.
     \end{displaymath}
     Here the ``constant'' by which the families fail to satisfy the
     condition depends on the multiplicity $m$, so that with respect
     to $m$ the asymptotics are not proper. However, we should like to
     point out that it does not
     depend on the number $r$ of singular points which are imposed.
     
     The families in Theorem \ref{thm:thm} thus exhibit the same
     properties as the families of plane curves provided in
     \cite{GLS01}, which we use to construct the non-T-smooth
     component. The idea is to intersect a family of cones in $\PC^3$
     over the plane curves provided by \cite{GLS01} with $\Sigma$ and
     to calculate the dimension of the resulting family. Under the
     conditions on $m$ and $l$ requested this family turns out to be
     of higher dimension than the expected one. The same idea was used
     by Chiantini and Ciliberto in \cite{CC99} in order to give
     examples of nodal families of curves on surfaces in $\PC^3$ which
     are not of the expected dimension. We then combine an asymptotic
     $h^1$-vanishing result by Alexander and Hirschowitz \cite{AH00}
     with an existence statement from \cite{KT02} to show that there
     is also a T-smooth component, where actually the curves have
     their singularities in very general position.

     \begin{proof}[Proof of Theorem \ref{thm:thm}]
       Fix a general plane $P$ in $\PC^3$ and a general point $p$.
       By \cite{GLS01} there is an integer $l_1=l_1(m)$ such that for
       any $l\geq \max\{l_1,m\}$ the
       family of curves in $P$ of degree $2lm+l$
       with $4l^2$ ordinary $m$-fold points as only singularities has
       a component $W$ of dimension 
       \begin{align*}
         \dim(W)\;\geq &\; (m+1)\cdot
         \frac{(l+1)\cdot(l+2)}{2}+(2l+1)\cdot(2l+2)-4\\
         =&\;\frac{l^2m+9l^2+3lm+15l+2m-4}{2}.
       \end{align*}
       Let $\kw$ be the family of cones with vertex $p$ over curves $C$ in $W$,
       then $\dim(\kw)=\dim(W)$, since a cone is uniquely determined by
       the curve $C$ and the vertex $p$. Moreover, any cone in $\kw$
       has precisely $4l^2$ lines of multiplicity $m$, so that when we intersect
       it with $\Sigma$ we get in general an irreducible curve in
       $\Sigma$ with $4l^2n$ ordinary $m$-fold points. In particular,
       $V_{|(2lm+l)H|}^{irr}\big(4l^2nM_m\big)$ must have a component $W'$ of dimension 
       \begin{displaymath}
         \dim(W')\geq\dim(W)\geq\frac{l^2m+9l^2+3lm+15l+2m-4}{2}.
       \end{displaymath}
       However, since the dimension of the linear system $|dH|$ is
       \begin{displaymath}
         \dim|dH|=\binom{d+3}{3}-\binom{d+3-n}{3}-1,
       \end{displaymath}
       and since 
       \begin{displaymath}
         \tau^{es}(M_m)=\frac{m\cdot(m+1)}{2}-2
       \end{displaymath}
       is the expected number of conditions imposed by an ordinary
       $m$-fold point, the expected dimension of
       $V_{|(2lm+l)H|}^{irr}\big(4l^2nM_m\big)$ is 
       \begin{align*}
         \expdim\Big(V_{|(2lm+l)H|}^{irr}\big(4l^2nM_m\big))\Big)= &
         \dim|(2lm+l)H|-4l^2n\tau^{es}(M_m)\\
         =&\frac{17l^2n+\big(4n-n^2\big)\cdot l\cdot (2m+1)}{2}+\frac{n^3-6n^2+11n}{6}.
       \end{align*}
       Due to the conditions on $m$ and $l$ this number is strictly
       smaller than the dimension of $W'$. It remains to show that
       $V_{|(2lm+l)H|}^{irr}\big(4l^2nM_m\big)$ also has a T-smooth component, after
       possibly enlarging $l_1$. 

       For 
       $\underline{z}=(z_1,\ldots,z_r)\in\Sigma^r$ we denote by
       $X(m;\underline{z})$ the zero-dimensional scheme with ideal
       sheaf $\kj_{X(m;\underline{z})}$ given by the stalks 
       \begin{displaymath}
         \kj_{X(m;\underline{z}),z}=
         \left\{
           \begin{array}{ll}
             \m_{\Sigma,z}^m,&\mbox{ if } z\in\{z_1,\ldots,z_r\},\\
             \ko_{\Sigma,z}, &\mbox{ else,}
           \end{array}
         \right.
       \end{displaymath}
       where $\ko_{\Sigma,z}$ denotes the local ring of $\Sigma$ at
       $z$ and $\m_{\Sigma,z}$ is its maximal ideal.

       By \cite{AH00} Theorem 1.1 there is an integer $l_2=l_2(m,\Sigma)$
       such that for $l\geq l_2$ and $\underline{z}\in\Sigma^r$ in
       very general position the canonical map
       \begin{displaymath}
         H^0\big(\Sigma,\ko_\Sigma((2lm+l-1)H)\big)\longrightarrow
         H^0\big(\Sigma,\ko_{X(m;\underline{z})}((2lm+l-1)H)\big)
       \end{displaymath}
       has maximal rank. In particular, since
       $h^1\big(\Sigma,\ko_\Sigma((2lm+l-1)H)\big)=0$ we have 
       \begin{displaymath}
         h^1\big(\Sigma,\kj_{X(m;\underline{z})}(2lm+l-1)H)\big)=0,
       \end{displaymath}
       once $\deg\big(X(m;\underline{z})\big)\leq
       h^0\big(\Sigma,\ko_\Sigma((2lm+l-1)H)\big)$, which is
       equivalent to
       \begin{displaymath}
         \frac{4l^2n\cdot m\cdot(m+1)}{2}\leq
         \binom{2lm+l+2}{3}-\binom{2lm+l+2-n}{3},
       \end{displaymath}
       or alternatively
       \begin{displaymath}
         \frac{nl\cdot\big(l-(n-2)\cdot (2m+1)\big)}{2}+\frac{n^3-3n^2+2n}{6}\geq 0.
       \end{displaymath}
       The latter inequality is fulfilled as soon as $l\geq
       (n-2)\cdot (2m+1)$. Moreover, under this hypothesis we
       have 
       \begin{displaymath}
         (2lm+l)\cdot H^2-2g(H)=
         (2lm+l)\cdot n-(n-1)\cdot(n-2)
         \geq 
         2m,
       \end{displaymath}
       where $g(H)$ denotes the geometric genus of $H$, and 
       \begin{equation}\label{eq:thm:1}
         (2lm+l)^2\cdot H^2>4l^2nm^2.
       \end{equation}
       Thus \cite{KT02} Theorem 3.3 (see also \cite{Kei01} Theorem
       1.2) implies that the family
       $V_{|(2lm+l)H|}^{irr}\big(4l^2nM_m\big)$ has a non-empty
       T-smooth component, more precisely it contains a curve in a
       T-smooth component with singularities in $z_1,\ldots,z_r$. In
       particular, since there is only a finite number of components
       and $\underline{z}$ is in very general position,
       some T-smooth component must dominate
       $\Sym^r(\Sigma)$. Actually, due to \cite{Los98} Proposition 2.1
       (e) and since $h^1(\Sigma,\ko_\Sigma)=0$ every T-smooth
       component dominates $\Sym^r(\Sigma)$.

       Thus the statement follows with
       \begin{displaymath}
         l_0(m,\Sigma):=\max\big\{l_1(m),l_2(m,\Sigma),
         (\deg(\Sigma)-2)\cdot (2m+1),m\big\}.
       \end{displaymath}

       It just remains to show that the fundamental group of the
       complement of a curve $C\in V_{|dH|}^{irr}(rM_m)$ is abelian. 
       Note first of all that by the Lefschetz Hyperplane Section
       Theorem $\Sigma$ is simply connected. But then
       $\pi_1(\Sigma\setminus C)$ is abelian by \cite{Nor83} 
       Proposition~6.5 because of \eqref{eq:thm:1}.
     \end{proof}


\providecommand{\bysame}{\leavevmode\hbox to3em{\hrulefill}\thinspace}

\end{document}